\documentclass[11pt,letterpaper]{amsart}

\usepackage[no-math]{fontspec}
\usepackage{microtype,mathtools,amssymb,enumitem,tikz-cd}
\usepackage[colorlinks=true,linkcolor=blue,citecolor=blue,urlcolor=blue]{hyperref}

\theoremstyle{plain}
\newtheorem{thm}{Theorem}[section]
\newtheorem{prop}[thm]{Proposition}
\newtheorem{lemma}[thm]{Lemma}

\theoremstyle{definition}
\newtheorem{rmk}[thm]{Remark}

\numberwithin{equation}{section}

\newcommand{\Pp}{\mathbb P}
\newcommand{\A}{\mathbb A}
\newcommand{\Ff}{\mathbb F}
\newcommand{\Oo}{\mathcal O}
\newcommand{\Q}{\mathbb Q}

\DeclareMathOperator{\Br}{Br}

\title[Differential forms and approximation on a K3 surface]
{The Brauer--Manin obstruction to weak approximation on a constant K3 surface}

\author{Domenico Valloni}
\date{}

\begin{document}

\begin{abstract}
We construct a constant supersingular K3 surface over every global field
of characteristic two for which rational points are Zariski dense but
are not dense in the Brauer--Manin set. This shows that the analogue
of Skorobogatov's conjecture over global function fields does not hold
for supersingular K3 surfaces.
\end{abstract}

\maketitle

\section{Introduction}

Rational points on K3 surfaces remain poorly understood.
A central conjecture of Skorobogatov \cite{Skorobogatov} asserts that the
Brauer--Manin obstruction is the only obstruction to the Hasse principle
and weak approximation for K3 surfaces over number fields.
In its weak approximation form, it predicts
$\overline{X(K)}=X(\A_K)^{\Br(X)}$.
The same statement can be formulated over global function fields.

In positive characteristic, supersingular K3 surfaces behave differently
from those of finite height. For K3 surfaces over number fields, and for
finite-height K3 surfaces over global function fields, the Brauer group
modulo constant classes is finite \cite{SZ,LS}.
For a supersingular K3 surface, the geometric Brauer group is instead
isomorphic to the additive group $\mathbb G_a(\overline K)$.
Over a global function field, the quotient by constant classes becomes
not finitely generated, after a suitable finite extension
\cite[Theorem~A and Remark~1.4]{LS}. The Brauer--Manin obstruction for supersingular K3 surfaces is studied
in \cite[Theorem~1.2]{Valloni} using differential forms. This result shows that,
after a suitable purely inseparable extension of the ground field, the
local projections of the Brauer--Manin set have measure zero at every
place. This is considerably stronger than in the number-field and
finite-height cases, where the Brauer--Manin set is open.

In this article, we prove that the
Brauer--Manin obstruction is not the only obstruction to weak approximation
for a supersingular K3 surface of Artin invariant one in characteristic
two. The example is constant, works over every global field of
characteristic two, and has Zariski dense rational points.
More precisely, let $k$ be a finite field of characteristic two and
consider the surface
\begin{equation}\label{surface}
 X=\left\{\sum_{i=0}^2x_i y_i^2=0,\quad
              \sum_{i=0}^2x_i^2 y_i=0\right\}
 \subset\Pp^2_k\times\Pp^2_k.
\end{equation}
This is the supersingular K3 surface of Artin invariant one studied by
Dolgachev and Kond\={o} \cite{DK}. Let $B/k$ be a smooth projective
geometrically integral curve and put $K=k(B)$. We write $X_K$ for the
base change of $X$ to $K$. For a field $L$ of characteristic two, write
$\Pp^2(L)^2=\{(u_0^2:u_1^2:u_2^2):(u_0:u_1:u_2)\in\Pp^2(L)\}$.
For each place $v$ of $K$, put
\[
 \begin{aligned}
 S_{x,v}&=\{(x,y)\in X(K_v):x\in\Pp^2(K_v)^2\},\\
 S_{y,v}&=\{(x,y)\in X(K_v):y\in\Pp^2(K_v)^2\}.
 \end{aligned}
\]
Define $S_x=\prod_v S_{x,v}$, $S_y=\prod_v S_{y,v}$, and
$D=\prod_v(S_{x,v}\cup S_{y,v})$.

\begin{thm}\label{main}
The set $X(K)$ is Zariski dense, and
\[
 D\subseteq X(\A_K)^{\Br(X_K)},\qquad
 \overline{X(K)}\cap D=S_x\cup S_y\subsetneq D.
\]
In particular, $\overline{X(K)}\subsetneq X(\A_K)^{\Br(X_K)}$.
If $K=k(t)$, then
\[
 X(\A_K)^{\Br(X_K)}=D,
 \qquad
 \overline{X(K)}=S_x\cup S_y.
\]
\end{thm}

\begin{rmk}\label{nowhere dense}
When $K=k(t)$, the closure
$\overline{X(K)}$ is nowhere dense in $X(\A_K)^{\Br(X_K)}$.
Its complement is therefore open and dense in the Brauer--Manin set.
\end{rmk}

The idea of the proof consists in geometrizing the differential
Brauer--Manin condition from \cite[Theorem~1.1]{Valloni}.
Since $X$ has Artin invariant one, the computations of van der Geer
and Katsura \cite[Theorem~6.2 and Corollary~11.12]{GK} show that
$H^0(X,\Omega^1_X/B_{1,X})$ is one-dimensional over $k$;
see also \cite[\S2.5.2]{Valloni}.
We construct an explicit nonzero class
$\omega\in H^0(X,\Omega^1_X/B_{1,X})$.
By \cite[Proposition~2.9]{Valloni}, there is a natural homomorphism $H^0(X_K,\Omega^1_{X_K}/B_{1,X_K})\longrightarrow\Br(X_K)[2].$
For the present surface, this map is surjective after base change
to an algebraic closure, by
\cite[Corollary~2.21 and \S2.5.2]{Valloni}.
For $c\in K$, multiplication by
$c^2$ preserves $B_{1,X_K}$, so $c^2\omega$ defines a section of this
quotient on $X_K$; we denote its image under this homomorphism by
$[c^2\omega]\in\Br(X_K)$.
Then, an adelic point
$(x_v)\in X(\A_K)$ is orthogonal to all the classes
$[c^2\omega]$, $c\in K$, if and only if the family
$(C(\omega_{|x_v}))_v$ comes from a global differential in $\Omega^1_K$.
Since $X$ and $\omega$ are defined over $k$, the argument in the
proof of \cite[Lemma~3.2]{Valloni}, applied to the constant model
$X\times_k B$, shows that this global differential is regular on $B$.
Let us assume that $K=k(t)$ for the sake of the introduction.
In this case $(C(\omega_{|x_v}))_v$ must vanish, and
orthogonality is equivalent to $C(\omega_{|x_v})=0$ at every place.
If $L/k$ is a field such that $[L:L^2]=2$, we show that an $L$-point
$a=(a_x,a_y)\in X(L)$ satisfies
$\omega_{|a}=0\in\Omega^1_L/B_{1,L}$ if and only if either
$a_x\in\Pp^2(L)^2$ or $a_y\in\Pp^2(L)^2$
(Lemma~\ref{form}). Adjoining square roots of the coordinates then
gives two purely inseparable covers $f_x:R_x\to X$ and $f_y:R_y\to X$.
Thus the two covers arise directly from the differential
Brauer--Manin condition. Their rationality will allow us to pass
from this differential description to a statement about the
full Brauer--Manin set.

\begin{rmk}\label{Hilbert property}
When $K=k(t)$, the discussion also gives
\[
 f_x(R_x(K))\cup f_y(R_y(K))=X(K).
\]
Over a number field, the Hilbert property precludes covering the
rational points by the images of finitely many finite covers of
degree greater than one. The covers here are purely inseparable.
Corvaja and Zannier \cite{CZ} proved the Hilbert property for a
diagonal quartic K3 surface over $\Q$. Gvirtz-Chen and Mezzedimi
\cite[Theorem~1.1 and \S1]{GM} proved that it holds, after a finite
extension, for every K3 surface of geometric Picard rank at least
$13$ over a number field, using the method developed by Demeio
\cite{Demeio}.
\end{rmk}

We check that both $R_x$ and $R_y$ are rational surfaces, and that
both covers have $k$-points above the same point of $X$.
Their rationality therefore shows that every Brauer class pulls back
to the same constant class on the two covers. Thus, when
$K=k(t)$, an adelic point belongs to $X(\A_K)^{\Br(X_K)}$ if and only
if each component lifts to a $K_v$-point of $R_x$ or $R_y$.
Over $k(t)$, a rational point must lift to one cover globally, whereas
the Brauer--Manin condition allows the choice to vary from place to place.
This explains the additional obstruction to weak approximation.
For a general base curve, we use the finiteness of
$H^0(B,\Omega^1_{B/k})$ to show that every point of
$\overline{X(K)}\cap D$ belongs to $S_x\cup S_y$.

\begin{rmk}
When $K=k(t)$, points of $X(K)$ correspond to morphisms
$g:\Pp^1_k\to X$, so our argument studies the constraints on rational
curves imposed by the global class $\omega$.
This follows the classical strategy of using global differential
tensors to constrain rational curves, as in Bogomolov's work
\cite{Bogomolov}. The feature specific to positive characteristic
that we exploit is that the sheaves $\Omega^1_X/B_{n,X}$ may have
nonzero global sections even when $H^0(X,\Omega^1_X)=0$, as it happens
for supersingular K3 surfaces.
\end{rmk}

\subsection{Some literature over number fields}

We briefly mention a few results related to the questions considered
here. This account is selective and incomplete. Over number fields, Skorobogatov and Swinnerton-Dyer \cite{SSD}
and Harpaz and Skorobogatov \cite{HS} obtained conditional results
towards the Hasse principle for certain families of Kummer surfaces
and varieties. Their results assume the finiteness of the relevant
Tate--Shafarevich groups, together with additional hypotheses on the
underlying abelian varieties. Results on $p$-adic density were obtained by Swinnerton-Dyer
\cite{SwinnertonDyer} and Pannekoek \cite{Pannekoek}.

\subsection{Notations}

We use absolute K\"ahler differentials and write $B_{1,T}=d\Oo_T$.
Over a perfect ground field, these agree with the relative differentials.
Evaluation of a class $\omega\in H^0(X,\Omega^1_X/B_{1,X})$ at an
$L$-point is denoted by $\omega_{|P}\in\Omega^1_L/B_{1,L}$.
The Cartier operator is denoted by $C$; for its properties and the
higher exact differential sheaves $B_{n,T}$, we follow
\cite[\S2.1]{Valloni}. All closures are taken in the adelic topology.

\section{A differential form and two rational covers}

The Jacobian criterion shows that \eqref{surface} is smooth: the first
equation has nonzero differential in the $x$-directions, and the second
has nonzero differential in the $y$-directions. Adjunction gives
$\Omega^2_X\simeq\Oo_X$, and the Koszul resolution gives
$H^1(X,\Oo_X)=0$. Thus $X$ is a K3 surface.

On the affine chart $x_0y_2\neq0$, write
$x=(1,a,b)$ and $y=(a^2c+b^2,c,1)$, using the second equation
to eliminate $y_0$. The first equation then gives
\begin{equation}\label{affine equations}
 b^4+b=(a^4+a)c^2,\qquad db=c^2\,da.
\end{equation}
In particular, $a,c$ are \'{e}tale coordinates on this chart.

\begin{lemma}\label{form}
There is a nonzero class $\omega\in H^0(X,\Omega^1_X/B_{1,X})$
represented by $a\,dc$ on this chart. If $L/k$ satisfies $[L:L^2]=2$,
then for every $P=(x,y)\in X(L)$,
\[
 \omega_{|P}=0
 \quad\Longleftrightarrow\quad
 x\in\Pp^2(L)^2\quad\text{or}\quad y\in\Pp^2(L)^2.
\]
\end{lemma}

\begin{proof}
The six opens $U_{ij}=\{x_i y_j\neq0\}$, with $i\neq j$, cover $X$.
On $U_{ij}$, let $h$ be the remaining index and take the class of
$(x_h/x_i)\,d(y_h/y_j)$.
Expressed in the function field using \eqref{affine equations}, these
representatives are $a\,dc$, $(b/c^2)\,dc$, or
$(ab/(a^2c+b^2)^2)\,dc$.
Their differences from $a\,dc$ are
\[
 \left(\frac{a^2c+b^2}{c}\right)^2dc
 \quad\text{or}\quad
 \left(\frac{ac}{a^2c+b^2}\right)^2dc,
\]
both exact. Since $\Omega^1_X/B_{1,X}$ is locally free over $\Oo_X^2$,
the classes agree on the overlaps and glue. Moreover,
$d\omega=da\wedge dc\neq0$.
Choose $t\in L$ with $L=L^2\oplus tL^2$ and write uniquely
$a=A_0^2+tA_1^2$ and $c=C_0^2+tC_1^2$. Then
\begin{equation}\label{Cartier calculation}
 C(a\,dc)=A_1C_1\,dt.
\end{equation}
Thus $\omega_{|P}=0$ if and only if $a\in L^2$ or $c\in L^2$.
If $a\in L^2$, then $db=c^2\,da=0$, hence $b\in L^2$,
so $x\in\Pp^2(L)^2$. If $c\in L^2$, then
$a^2c+b^2\in L^2$, so $y\in\Pp^2(L)^2$.
The converse implications are immediate from the coordinates.
The calculation on the other charts is identical.
\end{proof}

We now geometrize the two alternatives in Lemma~\ref{form}.
For the first, set $x=u^2$ and take the reduced equations.
The first equation of $X$ becomes
$(\sum_i u_i y_i)^2=0$, and the second becomes
$\sum_i u_i^4y_i=0$. This gives
\begin{equation}\label{rational cover}
 R_x=\left\{\sum_{i=0}^2u_i y_i=0,
                  \sum_{i=0}^2u_i^4 y_i=0\right\}
 \subset\Pp^2_u\times\Pp^2_y,
 \qquad f_x(u,y)=(u^2,y).
\end{equation}
Here $u^2=(u_0^2:u_1^2:u_2^2)$. Interchanging the two factors gives
$f_y:R_y\to X$.

\begin{lemma}\label{covers}
The surfaces $R_x$ and $R_y$ are smooth, projective, and $k$-rational.
The maps $f_x$ and $f_y$ are finite purely inseparable morphisms of degree
two. For every field extension $L/k$, their images on $L$-points are
$\{(x,y)\in X(L):x\in\Pp^2(L)^2\}$ and
$\{(x,y)\in X(L):y\in\Pp^2(L)^2\}$, respectively.
\end{lemma}

\begin{proof}
The projection $R_x\to\Pp^2_u$ determines $y$ uniquely unless the two
linear equations in $y$ are dependent. On a chart $u=(1,r,s)$, this
occurs exactly when $r^4=r$ and $s^4=s$. Thus the exceptional locus
in $\Pp^2_u$ is the reduced set of $21$ geometric points
$\Pp^2(\Ff_4)$. Near such a point, eliminating $y_0$ leaves the equation $(r^4+r)y_1+(s^4+s)y_2=0.$
Since $(r,s)\mapsto(r^4+r,s^4+s)$ is \'{e}tale, this is the usual
blow-up chart at a point. Hence $R_x$ is the blow-up of $\Pp^2$ at that
reduced finite subscheme. In particular, it is smooth and rational. The map $f_x$ is finite, being a closed restriction of Frobenius on the
first factor. On the chart of \eqref{affine equations}, if $r^2=a$, then
$s=b^2+(a^2+r)c$ satisfies $s^2=b$. Thus the function-field extension is
generated by a square root of $a$; since $da\neq0$, its degree is two.
For the assertion about points, let $P=(x,y)\in X(L)$ with $x=u^2$
for some $u\in\Pp^2(L)$. The first equation of $X$ gives
$(\sum_i u_i y_i)^2=0$, and hence $\sum_i u_i y_i=0$.
The second equation gives $\sum_i u_i^4y_i=0$.
Thus $(u,y)\in R_x(L)$ and $f_x(u,y)=P$.
Conversely, every point in $f_x(R_x(L))$ has square $x$-coordinates
by definition of $f_x$. The argument for $R_y$ is the same.
\end{proof}

\section{The adelic sets}

Consider the $k$-point $P_*=((1:0:0),(0:0:1)) \in X(k)$.
It has a $k$-point above it on each $R_x$ and $R_y$.
Put
\[
 B_\omega=\langle[c^2\omega]:c\in K\rangle
 \subseteq\Br(X_K)[2].
\]
The differential criterion controls orthogonality to $B_\omega$,
but we do not need to know whether these classes generate
$\Br(X_K)$ modulo constants. We first show directly that every
point of $D$ is orthogonal to the full Brauer group. This is
where we use that both surfaces $R_x$ and $R_y$ are rational.

\begin{lemma}\label{Brauer inclusion}
For every $\alpha\in\Br(X_K)$, every place $v$, and every
\[
 P_v\in S_{x,v}\cup S_{y,v}
 =f_x(R_x(K_v))\cup f_y(R_y(K_v)),
\]
one has
$\alpha_{|P_v}=\operatorname{res}_{K_v/K}(\alpha_{|P_*})\in\Br(K_v)$,
where $\operatorname{res}_{K_v/K}:\Br(K)\to\Br(K_v)$ is the natural map. In particular,
$D\subseteq X(\A_K)^{\Br(X_K)}$.
\end{lemma}

\begin{proof}
Since $R_x$ is smooth, projective, and rational,
$\Br((R_x)_K)=\Br(K)$; see \cite{CTS} for birational invariance of
the Brauer group of smooth proper varieties.
Evaluating at the point above $P_*$ shows that $f_x^*\alpha$
is the constant class $\alpha_{|P_*} \in \Br(K)$. The same holds for $f_y$.
Lemma~\ref{covers} gives the asserted equality of evaluations.
For an adelic point $(P_v)\in D$, global reciprocity then gives
\[
 \sum_v\operatorname{inv}_v\bigl(\alpha(P_v)\bigr)
 =
 \sum_v\operatorname{inv}_v
 \bigl(\operatorname{res}_{K_v/K}\alpha(P_*)\bigr)=0.
\]
Thus $D\subseteq X(\A_K)^{\Br(X_K)}$.
\end{proof}

For $P\in X(K)$, properness extends $P$ to a morphism $g_P:B\to X$.
Cartier therefore gives a regular differential
\begin{equation}\label{global differential}
 \eta_P=C(g_P^*\omega)\in H^0(B,\Omega^1_{B/k}).
\end{equation}
Lemma~\ref{form}, applied to $K$, gives
\begin{equation}\label{global zero}
 \eta_P=0
 \quad\Longleftrightarrow\quad
 P\in f_x(R_x(K))\cup f_y(R_y(K)).
\end{equation}
More generally, the differential criterion of
\cite[Theorem~1.1]{Valloni} associates to each
$(P_v)\in X(\A_K)^{B_\omega}$ a unique global differential $\eta$
whose image in $\Omega^1_{K_v}$ is $C(\omega_{|P_v})$ for every $v$.
As explained in the introduction, the constant model $X\times_k B$
shows that $\eta$ is regular on $B$. In particular, restriction
to the full Brauer--Manin set defines a map
\[
 \psi\colon X(\A_K)^{\Br(X_K)}\longrightarrow
 H^0(B,\Omega^1_{B/k})
\]
The target is finite because $k$ is finite. Lemma~\ref{form}
identifies vanishing of all the local evaluations with membership
in $D$, and Lemma~\ref{Brauer inclusion} shows that every point of
$D$ belongs to the domain of $\psi$. Hence
\begin{equation}\label{zero fibre}
 D=\psi^{-1}(0).
\end{equation}

\begin{prop}\label{closure}
For every $K=k(B)$, one has
$\overline{X(K)}\cap D=S_x\cup S_y$.
\end{prop}

\begin{proof}
The rational surfaces $R_x$ and $R_y$ satisfy weak approximation.
Their adelic spaces are compact, so their images $S_x$ and $S_y$ are
closed. Consequently,
\[
 \overline{f_x(R_x(K))}=S_x,\qquad
 \overline{f_y(R_y(K))}=S_y.
\]
This proves one inclusion. Let now $Q\in\overline{X(K)}\cap D$ and fix a place $v$.
The map $P_v\mapsto C(\omega_{|P_v})$ is continuous by
\cite[Proposition~3.10]{Valloni} and the continuity of the Cartier
operator. It vanishes at $Q_v$. The space $H^0(B,\Omega^1_{B/k})$
is finite as a set, and its
map to $\Omega^1_{K_v}$ is injective. Choose an open neighbourhood
$W$ of zero in $\Omega^1_{K_v}$ containing no image of a nonzero
element of this finite set, and let $U_v$ be its inverse image under
$P_v\mapsto C(\omega_{|P_v})$. Then $Q_v\in U_v$, and every global
point $P$ with $P_v\in U_v$ satisfies $\eta_P=0$.
Intersecting any adelic neighbourhood of $Q$ with the condition
$P_v\in U_v$ still gives a neighbourhood of $Q$, so it contains
a global point with $\eta_P=0$. By \eqref{global zero},
$Q$ therefore lies in the closure of
$f_x(R_x(K))\cup f_y(R_y(K))$, which is $S_x\cup S_y$.
\end{proof}

To prove that $S_x\cup S_y\subsetneq D$, we now construct an adelic
point in $D\setminus(S_x\cup S_y)$. There are many ways to do this. For example, the formula
\begin{equation}\label{curves}
 c_x(s)=\bigl((1:s^2:s^4),(s^5:s^4+s:1)\bigr)
\end{equation}
defines a morphism $\Pp^1\to X$; let $c_y$ be obtained by interchanging
the two factors. Notice that $c_x(0)=c_y(\infty)=P_*$.
Choose distinct places $v_0,v_1$, with uniformizers $\pi_0,\pi_1$, and set
\[
 Q_{v_0}=c_x(\pi_0),\qquad Q_{v_1}=c_y(\pi_1),\qquad
 Q_v=P_*\quad(v\neq v_0,v_1).
\]
At $v_0$, the $x$-projection lies in $\Pp^2(K_{v_0})^2$, whereas
$y_1/y_2=\pi_0^4+\pi_0$ has valuation one and is not a square.
At $v_1$, the reverse holds. Thus $Q\in D\setminus(S_x\cup S_y)$.
Together with Lemma~\ref{Brauer inclusion} and Proposition~\ref{closure},
this proves the strict inclusion in Theorem~\ref{main}. The failure of approximation at $Q$ can already be detected at
these two places. By the argument in the proof of
Proposition~\ref{closure}, there is an open neighbourhood $U_0$ of
$Q_{v_0}$ such that any global point $P$ with $P_{v_0}\in U_0$
satisfies $\eta_P=0$. The local cover images are closed, since the
local point spaces of $R_x$ and $R_y$ are compact. We may therefore
shrink $U_0$ to be disjoint from $S_{y,v_0}$ and choose an open
neighbourhood $U_1$ of $Q_{v_1}$ disjoint from $S_{x,v_1}$.
A global point meeting both conditions would, by \eqref{global zero},
lift to $R_x$ or $R_y$ over $K$. The condition at $v_0$ excludes an
$R_y$-lift, while that at $v_1$ excludes an $R_x$-lift.
Thus the adelic open set defined by $P_{v_0}\in U_0$ and
$P_{v_1}\in U_1$ contains the Brauer--Manin point $Q$ but no global
point. This explains why two places suffice. Finally, $R_x(K)$ is Zariski dense, so the dominant map $f_x$ shows
that $X(K)$ is Zariski dense as well.

If $K=k(t)$, then $B=\Pp^1_k$ and $H^0(B,\Omega^1_{B/k})=0$.
Thus \eqref{zero fibre} gives $X(\A_K)^{\Br(X_K)}=D$.
Also, $\eta_P=0$ for every global point, so \eqref{global zero}
and weak approximation on the covers give
$\overline{X(K)}=S_x\cup S_y$.
This completes the proof of Theorem~\ref{main}.

\begin{proof}[Proof of Remark~\ref{nowhere dense}]
Assume $K=k(t)$. By Theorem~\ref{main},
\[
 X(\A_K)^{\Br(X_K)}=D,\qquad
 \overline{X(K)}=S_x\cup S_y.
\]
Since $S_x\cup S_y$ is closed, it remains to show that it has
empty interior in $D$. Write $D_v=S_{x,v}\cup S_{y,v}$.
Since $X$ is proper, $D=\prod_v D_v$ carries the product topology.
Let $\mathcal U=\prod_v U_v$ be a nonempty basic open subset of $D$,
where each $U_v$ is open in $D_v$ and
$U_v=D_v$ outside a finite set $T$ of places.
Choose $A=(A_v)\in\mathcal U$ and distinct places $v_0,v_1\notin T$,
with uniformizers $\pi_0,\pi_1$. Replace only these two components
of $A$ by
\[
 Q_{v_0}=c_x(\pi_0)\in S_{x,v_0}\setminus S_{y,v_0},\qquad
 Q_{v_1}=c_y(\pi_1)\in S_{y,v_1}\setminus S_{x,v_1},
\]
and put $Q_v=A_v$ at every other place. Since $\mathcal U$ imposes
no restriction at $v_0$ or $v_1$, we still have $Q\in\mathcal U$.
However, its $v_0$-component excludes $Q\in S_y$, and its
$v_1$-component excludes $Q\in S_x$. Thus every nonempty basic
open subset of $D$ meets $D\setminus(S_x\cup S_y)$.
It follows that $S_x\cup S_y$ is closed with empty interior in $D$,
hence nowhere dense. By the two identifications above,
$\overline{X(K)}$ is therefore nowhere dense in
$X(\A_K)^{\Br(X_K)}$, as asserted in Remark~\ref{nowhere dense}.
\end{proof}

\begin{rmk}
For a base curve of positive genus, a global point need not belong to
either cover image. For example, let $B$ be the smooth projective curve
$b^2+b=a^3$. Then
$((1:a:b),(b:a:1))\in X(k(B))$ has neither projection in
$\Pp^2(k(B))^2$, and its differential \eqref{global differential} is
$C(a\,da)=da\neq0$. In general, the differential Brauer--Manin
criterion of \cite[Theorem~1.1]{Valloni} decomposes
$X(\A_K)^{B_\omega}$ into finitely many fibres indexed by
$H^0(B,\Omega^1_{B/k})$. The set $D$ is its zero fibre.
\end{rmk}
\section*{Final comments}

The role of small characteristic can also be seen in the geometry
behind the differential Brauer--Manin condition. For a superspecial
K3 surface $X$ over a field $k$ of characteristic $p$, with a nonzero
class $\omega\in H^0(X,\Omega^1_X/B_{1,X})$, the aim is to construct
a morphism $\pi\colon Y\to X$ such that, for field extensions $L/k$
with $[L:L^p]=p$,
\[
 \pi(Y(L))=
 \{P\in X(L):\omega_{|P}=0
       \text{ in }\Omega^1_L/B_{1,L}\}.
\]
A natural construction uses curvilinear subschemes of length $p$
in Frobenius fibres, subject to an exactness condition, and leads
to a parameter space $Y$ of dimension $p$. One can check that this
parameter space is projective for $p=2,3$, whereas for $p\geq5$
the curvilinear locus is nonproper. For the characteristic-two surface considered in this article,
one can take
\[
 Y=R_x\amalg R_y,\qquad \pi=f_x\amalg f_y.
\]
Thus the differential condition is expressed by two rational
surfaces. Their rationality is what makes the approximation
argument work. In characteristic three, one can check that the analogous
construction gives a smooth projective Calabi--Yau threefold
with $\ell$-adic Betti numbers
\[
 b_2=b_4=23,\qquad b_1=b_3=b_5=0
 \quad(\ell\ne3).
\]
Using the Kummer model over an algebraic closure, one can also
check that this threefold is birational to Schr\"oer's
characteristic-three Calabi--Yau threefold \cite{Schroer}.
The approximation argument used here does not immediately apply
to this threefold. In particular, we do not know whether the
Brauer--Manin obstruction is the only obstruction to weak
approximation on $Y$.

\section*{Use of artificial intelligence}

The author had been trying for some time to
`geometrise the Brauer--Manin obstruction' on supersingular K3
surfaces in every positive characteristic, using the methods of
\cite{Valloni}, without managing to obtain a clean statement.
In discussions of this problem, ChatGPT (OpenAI) suggested working
with the explicit Dolgachev--Kond\={o} model and, more importantly,
in characteristic two. The author had not realized the key role
played by small characteristic in the argument.

\end{document}